\documentclass[preprint,12pt]{elsarticle}
\usepackage{amsfonts}
\usepackage{amsmath}
\usepackage{amssymb}
\usepackage{amsthm}
\usepackage{color}
\usepackage{graphicx}%
\providecommand{\U}[1]{\protect\rule{.1in}{.1in}}
\newtheorem{thm}{Theorem}
\newtheorem*{thmA}{Theorem 1}
\newtheorem*{thmB}{Theorem 2}
\newtheorem*{thm3}{Theorem 3}
\newtheorem*{problem}{Problem}
\newtheorem{lem}[thm]{Lemma}
\newtheorem*{corolario}{Corollary}
\newdefinition{rmk}{Remark}
\newproof{pf}{\bf{Proof}}
\newproof{pot}{Proof of Theorem \ref{thm2}}
\newcommand{\QED}{\hfill {\qed}}
\begin{document}

\begin{frontmatter}

\title{On $H^{\infty}$ on the complement of $C^{1+\alpha}$ curves}
%\tnotetext[t1]{This document is a collaborative effort.}
%\tnotetext[t2]{The second title footnote which is a longer
%longer than the first one and with an intention to fill
%in up more than one line while formatting.}
%\author[esc]{J.M. Enr\'{i}quez-Salamanca\corref{cor2}\fnref{fn2}}

\author[esc]{J.M. Enr\'{i}quez-Salamanca\corref{cor1}}
\ead{enrique.desalamanca@uca.es}
\author[uca]{Mar\'{i}a J. Gonz\'{a}lez\corref{cor2}\fnref{fn2}}
\ead{majose.gonzalez@uca.es}
\cortext[cor2]{Corresponding author}
\cortext[cor1]{Principal corresponding author}
\address[uca]{Department of Mathematics, University of C\'{a}diz, Puerto Real, 11510, Spain}
\address[esc]{Department of Mathematics, University of C\'{a}diz, C\'{a}diz, 11002, Spain}
\fntext[fn2]{This author is supported by Grant MTM 2011-24606.}
\begin{abstract}
Let $\rho$ be a quasiconformal mapping on the plane with complex dilatation $\mu$. We show that if $\mu$ satisfies a certain Carleson measure condition, then one can transfer $H^{\infty}$ on the upper half plane onto the corresponding space in the complement of the quasicircle $\Gamma=\rho(\mathbb{R})$, and that this condition on $\mu$ characterizes $C^{1+\alpha}$ curves.

\end{abstract}

\end{frontmatter}
%\begin{document}
%\maketitle
\noindent {\bf Introduction}\\

 Let $f$ be a quasiconformal self-mapping of the plane with complex dilatation $\mu$. Thus $f$ is a homeomorphism with locally integrable distributional derivatives verifying that $\bar{\partial} f- \mu \partial f=0$, $\mu \in L^{\infty}(\mathbb{C})$ and $\|\mu\|_{\infty}<1$. \\
 % {\color{red} We shall also assume throughout this paper that $\mu$ has compact support.}\\

 The images of the real line (or the unit circle $\partial \mathbb{D}$) under quasiconformal mappings of the plane are called quasicircles. In general, they are not rectifiable and they do not satisfy any regularity conditions such as local absolute continuity or differentiability a.e., even when $\|\mu\|_{\infty}$ is small. Then, understanding the properties of the geometry of quasicircles in terms of the complex dilatation becomes one of the main objectives of the quasiconformal analysis and, also, of this article.\\

 One can get some regularity by imposing some stronger smallness condition on the dilatation $\mu$. If $\mu$ were zero in a neighborhood of $\partial\mathbb{D}$, then the map would be smooth on $\partial\mathbb{D}$. So, if $\mu$ decays to zero in some sense as it approaches the unit circle then we should be able to get some good behaviour of the mapping on $\partial\mathbb{D}$.\\

One of the first results along these lines is due to Carleson \cite{carleson}. He showed that if $f$ is a quasiconformal self mapping of $\mathbb{C} \backslash \mathbb{D}$ such that
$$\int_{0}^{1} \displaystyle{\frac{M(t)^2}{t}}dt < \infty,$$
\noindent where $M(t)=\mbox{sup}\{|\mu(z)|:1<|z|<1+t\}$, then $f$ is absolutely continuous on the circle and $f' \in L^{2}_{loc}$. Becker \cite{becker} extended this result to the case where $f$ represents a conformal mapping of $\mathbb{D}$ that extends quasiconformally to the whole plane.\\

In the same context, Dyn´kin \cite{dynkin} proved the following stronger result. Let

$$Q(\mu)(z)=\left(\int_{1<|\zeta|<2}\displaystyle{\frac{|\mu(\zeta)|^2}{|\zeta-z|^2}}d\xi
d\eta\right)^{1/2}.$$

\noindent If there exists a constant $a>0$ so that $e^{aQ(\mu)^2} \in L^{1}(\partial \mathbb{D})$, then the curve $\Gamma=f(\partial \mathbb{D})$ is rectifiable.\\

 A particular type of rectifiable quasicircle are the chord-arc or Laurentiev curves. A Jordan curve, $\Gamma$, is chord-arc if satisfies $\Lambda(\Gamma(z_1,z_2)) \leq C |z_1-z_2|$ for some constant $C>0$, where $\Lambda(\Gamma(z_1,z_2))$ denotes the length of the \mbox{shortest} arc of the curve $\Gamma$ between $z_1,z_2 \in \Gamma$.\\

 A condition for chord-arc curves with small constant was given by Astala-Zinsmeister (\cite{astala}, Theorem 3), and requires the measure

$$d \tau(z)=\displaystyle{\frac{|\mu(z)|^2}{|z|-1}}dx dy$$

\noindent to be a Carleson measure in $\mathbb{C} \backslash \mathbb{D}$ with small norm, i.e., $\|\tau\|_C \leq \gamma_0$ for some $\gamma_0>0$. The proof relies on estimates for the Schwarzian derivative of $f$ in $\mathbb{D}$, showing that $\log f' \in \operatorname{BMOA}$ with small $\operatorname{BMO}$ constant. On the other hand Semmes (\cite{SEMMES2}, Theorem 0.1) and MacManus (\cite{macmanus}, Corollary 6.5) showed that the same result holds with no assumption on $f$ being conformal on $\mathbb{D}$. Semmes' proof is based on estimates for a certain perturbed Cauchy integral operator while the strategy in \cite{macmanus} is to estimate Haar type coefficients for $\log f'$.\\

 For arbitrary constants, the result is no longer true \cite{bishop}. In fact, if no res\-tric\-tion on the Carleson norm of $|\mu(z)|^2/(|z|-1)$ is imposed, the quasicircle $\Gamma$ might not be even rectifiable.\\

 Related to these results is the following one due to Pommerenke \cite{pommerenke}. A curve $\Gamma$ is an asymptotically smooth curve, that is $$\displaystyle{\frac{\Lambda(\Gamma(z_1,z_2))}{|z_1-z_1|}}\rightarrow
      1, \;\;\;  \text{as} \;\;\; |z_1-z_2|\rightarrow 0,$$
%\end{thm}
\noindent if and only if $\log f' \in \operatorname{VMOA}.$\\

 This paper presents a new characterization of smooth curves in terms of $\mu$. In particular, we are interested in $C^{1+\alpha}$  curves. We prove the following result:
\begin{thmA}
\label{theoremA}
Let $f$ denote a conformal map of $\mathbb{D}$ onto the inner domain of a Jordan curve $\Gamma$. Then $\Gamma$ is a $C^{1+\alpha}$ curve if and only if $f$ extends to a global quasiconformal map whose dilatation $\mu$ satisfies that $|\mu(z)|^2/(|z|-1)^{1+\varepsilon}$ is a Carleson measure relative to $\mathbb{T}$, where  $\varepsilon=\varepsilon(\alpha)$ and $\alpha=\alpha(\varepsilon,\|\mu\|_{\infty})$.
\end{thmA}

 The equivalent result holds if we consider unbounded $C^{1+\alpha}$ curves. In this case $\mu$ satisfies that $|\mu(z)|^2/|y|^{1+\varepsilon}$ is a Carleson measure relative to $\mathbb{R}$, where $y=\operatorname{Im}(z)$.\\

 In the second part of the paper we consider quasiconformal mappings of the plane $\rho:\mathbb{C} \rightarrow \mathbb{C}$ whose complex dilatation $\mu$ satisfies that for some $\varepsilon >0$, $|\mu|^2/|y|^{1+\varepsilon}$ is a Carleson measure relative to $\mathbb{R}$. We will show that under this condition on $\mu$ we can transfer $H^{\infty}$ on the half plane $\mathbb{R}^{2}_{+}$ onto the corresponding space in the complement of the quasicircle $\Gamma=\rho(\mathbb{R})$. More precisely, denoting by $C_{\Gamma}(g)$ the Cauchy integral of a function $g\in L^{\infty}(\Gamma)$, that is

\begin{equation}
   C_{\Gamma}(g)(z)=\displaystyle{\frac{1}{2 \pi i}}\int_{\Gamma} \displaystyle{\frac{g(\omega)}{\omega - z}} d \omega, \; \; z \notin \Gamma, \nonumber
\end{equation}

\noindent we prove the following result:

 \begin{thmB}
 \label{theoremB}
 Let $\rho$ be a quasiconformal map of the plane onto itself whose complex dilatation $\mu$ has compact support and satisfies that for some $\varepsilon> 0$ , \newline $|\mu|^2/|y|^{1+\varepsilon}$ is a Carleson measure relative to $\mathbb{R}$. Let $\Omega_{+}$ and $\Omega_{-}$ denote the two regions bounded by the quasicircle $\Gamma=\rho(\mathbb{R})$ and let $g \in L^{\infty}(\Gamma)$. Then $C_{\Gamma}(g) \in H^{\infty}(\Omega_{\pm})$ if and only if $C_{\mathbb{R}}(f) \in H^{\infty}(\mathbb{R}^2_{\pm})$ respectively, where \mbox{$f=g\circ \rho$.}
 \end{thmB}

 The proof is based on an argument by Semmes \cite{SEMMES}, where the idea is to transform a $\bar{\partial}$ problem relative to $\Gamma$ into a $\bar{\partial}-\mu \partial$ problem relative to $\mathbb{R}$ via a change of variables. As an immediate consequence of Theorem 1 and Theorem 2 we obtain the following corollary.

%We will apply the previous results to study the relation between the Cauchy integral $C_{\Gamma}(g)$ of a function $g \in L^{\infty}(\Gamma)$, that is
%
%\begin{equation}
%   C_{\Gamma}(g)(z)=\displaystyle{\frac{1}{2 \pi i}}\int_{\Gamma} \displaystyle{\frac{g(\omega)}{\omega - z}} d \omega, \; \; z \notin \Gamma, \nonumber
%\end{equation}

%\noindent and the Cauchy integral of its pullback on $\mathbb{R}$. So, as an immediate consequence of Theorem 1 and Theorem 2 we obtain the following corollary.

\begin{corolario}
Let $\Gamma$ be an unbounded $C^{1+\alpha}$ curve analytic at $\infty$, and let $\rho$ denote a conformal map of $\mathbb{R}^2_+$ onto any of the regions bounded by $\Gamma$. Then, given a function $g \in L^{\infty}(\Gamma)$, the Cauchy integral $C_{\Gamma}(g) \in L^{\infty}(\mathbb{C})$ if and only if $C_{\mathbb{R}}(f) \in L^{\infty}(\mathbb{C})$, where $f$ denotes the pullback of $g$ under the conformal mapping $\rho$.
\end{corolario}

%%Although this integral might well diverge for $g \in L^{\infty}(\Gamma)$, it could always be defined module constants as the integral $G(z)-G(\omega)$ does make sense for any $z, \omega \in \mathbb{C} \backslash \Gamma$. More precisely, we prove the following result:
%
%\begin{thmB}
%\label{theoremB}
%  Let $\Gamma$ be a $C^{1+\alpha}$ curve and let $g \in L^{\infty}(\Gamma)$. Then, the Cauchy integral $C_{\Gamma}(g) \in L^{\infty}(\mathbb{C})$ if and only if $C_{\mathbb{R}}(f) \in L^{\infty}(\mathbb{C})$, where $f$ denotes the pullback of $g$ under the conformal mapping.
%\end{thmB}

 %The proof is based on an argument by Semmes \cite{SEMMES}, in which he relates a similar problem to the study of a certain $\bar{\partial}$ problem.\\

 The paper is structured as follows: In section 1, we review  some definitions and basic facts, in particular the analytic characterization of $C^{1+\alpha}$ curves. The proof of Theorem 1 is presented in section 2 whereas section 3 is devoted to Theorem 2.\\

% The authors would like to thank the referee for the careful revision of this article and for the many helpful suggestions to improve its reading.

\section{Preliminaries.}

 Let us denote complex variables by $z=x+iy$ and $\zeta=\xi+i\eta$. We shall use the following notation throughout this article: $\text{Im}(z)=y$, $\mathbb{D}=\{z:|z|<1\}$, $\mathbb{T}=\partial \mathbb{D}$, $B_x(r)$ denotes the ball centered at $x$ and radius $r$, $|I|$ represents the length of any arc $I\subset\partial\mathbb{D}$, $\Gamma(z_1,z_2)$ is the shortest arc of the curve $\Gamma$ between $z_1,z_2 \in \Gamma$, and $\Lambda(\Gamma)$ the length of the curve $\Gamma$. Also, we shall write $\bar{\partial}=\partial/\partial \bar{z}=1/2(\partial_x+i\partial_y)$ and $\partial=\partial/\partial z=1/2(\partial_x-i\partial_y)$.\\

A positive measure $\lambda$ on $\mathbb{C}$ is called a Carleson measure relative to a given chord-arc curve $\Gamma$ if there exists a constant $C>0$ such that $\lambda(B_{z}(R))\leq C\;R$ for all $z \in \Gamma$ and $R>0$. The smallest such $C$ is the norm of $\lambda$, $\|\lambda\|_C$. If \begin{equation}
\lim_{r \rightarrow 0} \; \sup_{R<r} \displaystyle{\frac{\lambda(B_{z}(R))}{R}}=0, \nonumber
\end{equation}
\noindent  the measure is said to be a vanishing Carleson measure or that it satisfies the o$(1)$-Carleson condition.\\
%$\lambda(B_{z}(R))$ goes to 0 as $R$ goes to 0, then we say that $\lambda$ is a vanishing Carleson measure.}\\

%\noindent The interior and exterior Carlenson´s squares are denoted by:
%$$S_i(I)=\{z=re^{i\theta}:e^{i\theta}\in I,1-|I| \leq r \leq
%1\},$$
%
%\noindent and
%
%$$S_e(I)=\{z=re^{i\theta}:e^{i\theta}\in I,1 \leq r \leq 1+|I|\},$$
%
%\noindent respectively. We will also consider the relation $a
%\asymp b$ if $c|b| \leq |a| \leq C|b|$, where $c$ and $C$ are
%positive constants , not necessarily the same throughout a
%formula.\\
%
% Let $\nu$ be a positive Borel measure in $\mathbb{D}$. The measure $\nu$ is said to be a Carleson measure if there exists
%a constant $C>0$ so that, $ \nu(S_i(I))\leq C \cdot |I|$ for every arc $I \subset \mathbb{T}$. The least constant $C$ is called the Carleson constant
%or the Carleson norm of the measure $\nu$, $\|\nu\|_{C}$.\\

%\noindent Also, let us define the limit:
%
%
%
%\begin{equation}
%N_0(\nu)=\displaystyle{\lim_{\delta \rightarrow 0}} \sup_{|I| <  \delta}
%\displaystyle{\frac{\nu(S_i(I))}{|I|}}. \nonumber
%\end{equation}
%
%\noindent If $N_0(\nu)=0$, then the measure $\nu$ is said to
%satisfy a vanishing Carleson condition or a $o(1)$- Carleson
%condition.\\
%
%The notion of Carleson measures are closely related to BMO functions \cite{garnett} .
We will denote by $H^{p}(\mathbb{D})$, $0 < p < \infty$, the Hardy space of analytic functions on $\mathbb{D}$ such that
 $$\sup_{r} \displaystyle{\frac{1}{2 \pi}}\int|f(re^{i\theta})|^p d\theta=\|f\|^p_{H_p}<+\infty.$$

 \noindent If $p=\infty$, $f \in H^{\infty}(\mathbb{D})$ if $f(z)$ is a bounded analytic function on $\mathbb{D}$, $\|f\|_{\infty}=\text{sup}_{z\in \mathbb{D}}|f(z)|$.\\

 A function $f\in L^1(\mathbb{T})$ belongs to the space $\operatorname{BMO}(\mathbb{T})$ if there exists $A>0$ so that
     $$\sup_{I}\displaystyle{\frac{1}{|I|}}\int_{I}|f(\zeta)-a_I||d\zeta| \leq A,
            \; \text{with}  \; a_I=\displaystyle{\frac{1}{|I|}}\int_{I}f(\zeta)|d\zeta|$$

\noindent and where the supremum is taken over all arcs $I \subseteq
   \mathbb{T}$. The least possible $A$ in this inequality is called the
$\operatorname{BMO}$ norm of $f$, $\|f\|_{*}$.\\

The space $\operatorname{VMO}(\mathbb{T})$ is defined as follows:
  $$\operatorname{VMO}(\mathbb{T})=\{f \in \operatorname{BMO}(\mathbb{T}): \lim_{|I|\rightarrow 0}\displaystyle{\frac{1}{|I|}}\int_{I}|f(\zeta)-a_I||d\zeta|=0\}.$$

\noindent We can extend the definitions of $\operatorname{BMO}$ and $\operatorname{VMO}$ to any locally rectifiable curve by replacing intervals with arcs.\\

We say that $f \in \operatorname{BMOA}$ or $f \in \operatorname{VMOA}$ if $f \in H^1(\mathbb{D})$ and if the boundary values of $f$ on $\mathbb{T}$
belong to $\operatorname{\operatorname{BMO}}(\mathbb{T})$ or $\operatorname{VMO}(\mathbb{T})$ respectively. Recall that $\operatorname{\operatorname{BMOA}}$ is contained in the Bloch space $$\mathcal{B}=\{f \;\mbox{analytic in}\; \mathbb{D}: \|f\|_{\mathcal{B}}=\mbox{sup}(1-|z|^2)|f'(z)|< \infty \}$$ and $\operatorname{VMOA}$ is contained in $\mathcal{B}_0 =\{f \in \mathcal{B}: (1-|z|^2)|f'(z)|\rightarrow 0 \;\mbox{as}\; |z|\rightarrow 1-0\}$. \\

The notion of Carleson measures is closely related to $\operatorname{BMO}$ functions (see for example Chapter 6, Section 3 in \cite{garnett}). Let $F$ be an analytic function on $\mathbb{R}^2_{+}$ and $f=F|_{\mathbb{R}}$. Then $f \in \operatorname{BMO}(\mathbb{R})$ if and only if $|F'(z)|^2 |y|\; dx dy$ is a Carleson measure with respect to $\mathbb{R}$ (\cite{garnett}, p 262). \\

%Analogous definitions and results extend to analytic functions defined on the complement of chord-arc curves \cite{jerison}. In this case, $y$ should be replaced by the distance of the point to the curve, $\delta_{\Gamma}(z)$.\\

%\begin{thm}
%\label{garnett1}
%  Let $f$  be an analytic function in $\mathbb{D}$. If the measure
%
%      $$d \nu=(1-|z|)|f'(z)|^2dxdy$$
%
%\noindent satisfies the Carleson condition in $\mathbb{D}$, then $f \in BMOA$. Furthermore, if $\nu$ satisfies the $o(1)-Carleson$ condition, then  $f(re^{i \theta})$ converges in $BMO$ norm and the boundary function belongs to $VMO$.
%
%\end{thm}

Let $\omega$ be a nonnegative locally integrable function on $\mathbb{R}$. Set $\omega(E)=\int_E \omega(x)dx$, and let $|E|$ denote the Lebesgue measure of $E$. We say that $\omega$ is an $A_\infty$ weight if for every $\varepsilon >0$, there exists $\delta >0$ such that if $I$ is any interval and $E \subseteq I$, then $|E|/|I| < \delta$ implies $\omega(E)/\omega(I) < \varepsilon$. Note that if $h$ is a bilipschitz map, clearly $|h'| \in A_{\infty}$. \\ %%Also, if $\varphi_1, \varphi_2$ are $A_\infty$ weights, then $\varphi_1 \circ \varphi_2$, $\varphi_1^{-1} \in A_{\infty}$.}\\

An important fact is that $\log \omega \in \operatorname{BMO}$ if $\omega \in A_{\infty}$, and $\{\log \omega:\omega \in A_{\infty}\}$ spans an open subset of real-valued $\operatorname{BMO}$, inducing a natural topology on $A_{\infty}$. In particular, there is a $\gamma >0$ so that $e^b \in A_{\infty}$ if $b$ is real valued and $\|b\|_{*} \leq \gamma$ (\cite{pommerenke}, p 171).\\

 We now introduce the curves which are the main object of study in this paper. A Jordan curve $\Gamma$ is said to be of class $C^{n}$ $(n=1,2,\ldots)$ if it has a parametrization $\varphi(\tau)=f(e^{i\tau}),0\leq \tau \leq 2\pi$, that is $n$ times continuously differentiable and satisfies that $\varphi'(\tau)\neq 0$, $\forall \tau$.  Furthermore, it is of class $C^{n+\alpha}$, for $0<\alpha <1$, if

 \begin{equation}
\label{curva}
 |\varphi^{(n)}(\tau_1)-\varphi^{(n)}(\tau_2)| \leq C|\tau_1-\tau_2|^{\alpha}.
 \end{equation}

It is well known that for $0 < \alpha < 1$ we can consider the parametrization of the curve given by the conformal mapping $f$ that sends $\mathbb{D}$ onto the inner domain bounded by $\Gamma$ (Kellogg-Warschawski Theorem). In this case, by the Hardy-Littlewood criterion (\cite{stein}, Section V.4), (\ref{curva}) is equivalent to

\begin{equation}
\label{pommerenke}
|f^{n+1}(z)| \leq C (1-|z|)^{\alpha-1}, \;\; \mbox{for all}\;\; z \in \mathbb{D}.
\end{equation}

%. A characterization of $C^{n,\alpha}$ curves was given by Pommerenke \cite{pommerenke}.
%\begin{thm}
% Let $f$ map $\mathbb{D}$ conformally onto the inner domain of
% a Jordan curve $\Gamma$ of class $C^{n,\alpha}$, where
% $n=1,2,\dots$ and $0<\alpha<1$. Then, $f^{(n)}$ has a
% continuous extension to $\bar{\mathbb{D}}$ and
%    \begin{equation}
%       |f^{(n)}(z_1)-f^{(n)}(z_2)|\leq C |z_1-z_2|^{\alpha}
%       \;\; \text{for} \; \; z_1,z_2 \in \bar{\mathbb{D}}.\nonumber
%    \end{equation}
%\end{thm}

%\noindent The conversal is trivial as we can always choose the
%conformal parametrization of the curve.\\
%\noindent METER LA OTRA CARACTERIZACIÓN DE CURVA $C^{1+\alpha}$ EN FUNCIÓN DE f´´.

\section{Proof of Theorem 1.}

The main idea is based on an estimate of the logarithmic derivative de\-ve\-lo\-ped by Dynkin \cite{dynkin}. Accordingly, if $f$ is a conformal mapping in the unit disc with a $k$-quasiconformal extension to the whole plane, that is $\|\mu\|_{\infty}\leq k < 1$, then for all $z \in \mathbb{D}$

\begin{equation}
\label{dynkin1}
(1-|z|)\left|\displaystyle{\frac{f''(z)}{f'(z)}}\right|\leq C(1-|z|)^{1-k}\left[1+\int_{1-|z|}^{1}\displaystyle{\frac{\omega(z,t)}{t^{2-k}}}dt\right],
\end{equation}
\noindent where
\begin{equation}
\label{omega}
\omega(z,t)=\left(\displaystyle{\frac{1}{\pi t^2}}\int_{|\zeta  -z|<t}|\mu(\zeta)|^2 d \xi d \eta\right)^{1/2},
\end{equation}

\noindent and where the constant $C$ depends on $k$ only.\\

\noindent Using this estimate, he also proved the following result (\cite{dynkin}, Theorem 2):

%\begin{equation}
%\int_{0}^{1}\displaystyle{\frac{\omega(z,t)}{t}}dt<\infty
%\end{equation}
%
%\noindent converges uniformly in $z \in \partial\mathbb{D}$, then $\log f'$ (and thereby $f'$ and $1/f'$) is continuous in the closed disc.
\begin{thm3}
\label{acotacion}
   If the integral

      \begin{equation}
        \int_{0}^{1}\displaystyle{\frac{\omega(z,t)}{t}}dt<\infty \nonumber
      \end{equation}

\noindent converges uniformly in $z \in \mathbb{T}$, then $log
f'$ (and thereby $f'$ and $1/f'$) is continuous in the closed disc.
\end{thm3}

\begin{thmA}
\label{theoremA}
Let $f$ denote a conformal map of $\mathbb{D}$ onto the inner domain of a Jordan curve $\Gamma$. Then $\Gamma$ is a $C^{1+\alpha}$ curve if and only if $f$ extends to a global quasiconformal map whose dilatation $\mu$ satisfies that $|\mu(z)|^2/(|z|-1)^{1+\varepsilon}$ is a Carleson measure relative to $\mathbb{T}$, where $\varepsilon=\varepsilon(\alpha)$ and $\alpha=\alpha(\varepsilon,\|\mu\|_{\infty})$.
\end{thmA}
%\begin{thmA}
%Let $f:\mathbb{D}\rightarrow \Omega$ be a Riemann map from $\mathbb{D}$ onto a region $\Omega$. $\Gamma=\partial \Omega$ is a $C^{1+\alpha}$ curve if and only if $f$ can be extended to a global quasiconformal map whose dilatation $|\mu|$ satisfies that $|\mu(z)|^2/(|z|-1)^{1+\varepsilon}$ is a Carleson measure, where $\varepsilon=\varepsilon(\|\mu\|_{\infty})$ and $\alpha=\alpha(\varepsilon,\|\mu\|_{\infty})$.
%\end{thmA}
\begin{pf}
\noindent Let us assume first that there exists a quasiconformal extension of $f$ so that for some $\varepsilon>0$, $\nu(z)=|\mu(z)|^2/(|z|-1)^{1+\varepsilon}$ is a Carleson measure relative to $\mathbb{T}$, $\|\mu\|_{\infty} \leq k <1$. Let $z\in \bar{\mathbb{D}}$, then for any $t>1-|z|$, we obtain from (\ref{omega}) that:

%By using this condition of Carleson measure, for any $t>1-|z|$, $z \in \mathbb{D}$, we obtain from (\ref{omega}) that:
\begin{eqnarray}
\label{acotacion7}
\omega^2(z,t)&
%=\displaystyle{\frac{1}{\pi t^2}}\int_{|\zeta  -z|<t}|\mu(\zeta)|^2 d \xi d \eta \leq &&\\
\leq & \displaystyle{\frac{1}{\pi
t^2}}\left(\int_{B_z(t)}\displaystyle{\frac{|\mu(\zeta)|^2}{(|\zeta|-1)^{1+\varepsilon}}}d
\xi d \eta\right)^{1/2}\left(\int_{B_{z}(t)} |\mu(\zeta)|^2
(|\zeta|-1)^{1+\varepsilon}d \xi d \eta\right)^{1/2}\nonumber\\
&\leq&  C t^{\varepsilon/2}.
\end{eqnarray}

%\begin{center}
%\begin{eqnarray}
%\omega(z,t)^2 & \leq & \displaystyle{\frac{1}{\pi t^2}}(\int_{|\zeta-z|<t}\displaystyle{\frac{|\mu(\zeta)|^2}{(1-|\zeta|)^{1+\varepsilon}}}d \xi d \eta)^{1/2}(\int_{|\zeta-z|<t}|\mu(\zeta)|^2 (1-|\zeta|)^{1+\varepsilon}d \xi d \eta)^{1/2} \leq \\
%              & \leq & \displaystyle{\frac{C}{t^{3/2}}}(\int_{|\zeta-z|<t}|\mu(\zeta)|^2 (1-|\zeta|)^{1+\varepsilon}d \xi d \eta)^{1/2} \leq
%              C(\varepsilon,k)t^{\varepsilon/2}.\\
%\end{eqnarray}
%\end{center}
\noindent where $C=C(\|\mu\|_{\infty},\|\nu\|_{C})$. To prove the last inequality note that the first integral can be always approximated by integrals on balls centered at the boundary. Then, by the Carleson condition on the measure $\nu$, we get that the first integral is bounded by $C(\|\nu\|_C) t^{1/2}$, while the second one is clearly bounded by $C(\|\mu\|_{\infty})t^{(3+\varepsilon)/2}$. Therefore, if $z \in \mathbb{D}$ and $\alpha \leq \min{(\varepsilon/4;1-k)}$, by (\ref{dynkin1})

   \begin{equation}
        \left|\displaystyle{\frac{f ''(z)}{f '(z)}}\right| \leq  C(1-|z|)^{-k}\left[1+C\int_{1-|z|}^{1}t^{\varepsilon/4+k-2}dt\right] \leq C(1-|z|)^{\alpha-1} \nonumber
   \end{equation}

\noindent and $|f ''(z)| \leq C|f '(z)| (1-|z|)^{\alpha-1}$, where now $C=C(\varepsilon,\|\mu\|_{\infty},\|\nu\|_{C})$.\\

%As $0\leq (1-|z|)\leq 1$ and $\alpha \leq 1-k$
%
%\begin{center}
%  \begin{equation}
%     \displaystyle{\frac{1}{(1-|z|)^{k}}}=\displaystyle{\frac{1}{(1-|z|)^{1-(1-k)}}} \leq \displaystyle{\frac{1}{(1-|z|)^{1-\alpha}}}
%  \end{equation}
%\end{center}

We need to prove that $|f'(z)|$ is
bounded on $\mathbb{T}$. Let $z \in \mathbb{T}$, then by (\ref{acotacion7})

  \begin{equation}
       \int_{0}^{1}\displaystyle{\frac{\omega(z,t)}{t}}dt \leq C\int_{0}^{1}\displaystyle{\frac{t^{\varepsilon/4}}{t}}dt< \infty,\nonumber
  \end{equation}

\noindent so we get by Theorem 3 that $|f'|$ is bounded in the closed disc. Therefore, $|f''(z)|\leq C(1-|z|)^{\alpha-1}$ which implies by  (\ref{pommerenke}) that $\Gamma$ is $C^{1+\alpha}$.\\

To prove the second part of the theorem, let us consider the following quasiconformal extension of the Riemann mapping f:

   \begin{equation}
    \label{extension}
       f(z)=f(1/\bar{z})+f'(1/\bar{z})(z-1/\bar{z}), \;\;\;\text{for}\;\;\;
       |z|>1.
    \end{equation}

Note that if $\beta(z)$ denotes the logarithmic derivative, that is $\beta(z)=(1-|z|)|f''(z)/f'(z)|$, $z \in \mathbb{D}$, and $\mu$ denotes the complex dilatation of the quasiconformal extension,  then for $|z|>1$

\begin{equation}
\label{acotacion5}
|\mu(z)| \asymp 1/|z|\beta(1/\bar{z}).
\end{equation}

In general, the mapping (\ref{extension}) is not homeomorphic, but Becker and Pommerenke (\cite{becker2}, Th. 4) proved that (\ref{extension}) is indeed a quasiconformal extension of f to a neighborhood of $\mathbb{T}$ if $f(\mathbb{D})$ is a Jordan domain and $\overline{\lim}_{|z|\rightarrow 1-0} \beta(z)<1$.\\

%\begin{equation}
%   f_z=f'(\displaystyle{\frac{1}{\bar{z}}})\;\;\;  \text{and} \;\;\;
%   f_{\bar{z}}=-\displaystyle{\frac{1}{\bar{z}^2}}(z-\displaystyle{\frac{1}{\bar{z}}})f''(\displaystyle{\frac{1}{\bar{z}}}). \nonumber
%   \end{equation}

% \noindent {\bf In general, this mapping is not homeomorphic,
%but Becker  \cite{becker3} proved that this
%extension was indeed a quasiconformal extension of $f$ to $\mathbb{C}$
%if $\overline{\lim}_{|z|\rightarrow 1-0}\beta(z)<1$.}\\

If $\Gamma$ is $C^{1+\alpha}$, it is assymptotically smooth, so by Pommerenke's result $\log f' \in \operatorname{VMOA}$ (see \cite{pommerenke}, pg. 172). Since $\operatorname{VMOA} \subset \mathcal{B}_0$, we get that
$\beta(z) \rightarrow 0$ as $|z| \rightarrow 1-0$ and the aforementioned extension of $f$ gives a well defined quasiconformal extension to a neighborhood of $\mathbb{T}$, $G=\{z:|z|<R_1 \text{for}\; \text{some} \;R_1>1\}$\\

To obtain a global quasiconformal mapping we apply the following theorem (\cite{lehto}, Th. 8.1): If $f:G\rightarrow G'$ is a $k$-qc map and $E$ is a compact set of the domain $G$, then there exists a $k'$-qc map of the whole plane that coincides with $f$ in $E$ and with $k'$ depending only on $k$, $G$ and $E$.\\

\noindent Setting $E=\{z:|z|\leq R_0, 1<R_0<R_1\}$, the above result provides a global quasiconformal extension of the conformal mapping $f$ whose complex dilatation $\mu$ satisfies (\ref{acotacion5}) for $1<|z|\leq R_0$.\\

\noindent Since $\Gamma$ is a $C^{1+\alpha}$ curve, $\log f'$ is continuous in $\bar{\mathbb{D}}$ (\cite{pommerenke}, Theorem 3.5) and therefore $|f'|$ is bounded below in $\mathbb{D}$. So, we get by (\ref{pommerenke}) and (\ref{acotacion5}) that

%\noindent It is enough to prove that the Carleson condition holds for balls $B_z(R_0)$ with $R_0$ small enough so that (\ref{acotacion5}) holds. As $|\mu(z)| \asymp 1/|z|^2 \beta(1/\bar{z})$ in a neighborhood of $\mathbb{T}$, $1<|z|\leq R_0$, and as $\Gamma$ is a $C^{1+\alpha}$ curve, by (\ref{pommerenke})
\begin{equation}
\label{acotacion6}
|\mu(z)|\leq C \beta(\displaystyle{\frac{1}{\bar{z}}}) \leq
C(|z|-1)^{\alpha}, \;\;\; \text{for} \; 1<|z|\leq R_0.
\end{equation}

To finish the proof of the theorem, it only remains to show the Carleson condition. For that, and without any loss of generality, consider a ball $B_z(R)$ centered at a point $z\in \mathbb{T}$ and radius $R<R_0-1$. By using (\ref{acotacion6}) and a change of variables to polar coordinates, we obtain that

$$\int_{B_z(R)}\displaystyle{\frac{|\mu(\zeta)|^2}{(|\zeta|-1)^{1+\varepsilon}}}d \xi d \eta \leq C \int_{B_z(R)}\displaystyle{\frac{(|\zeta|-1)^{2\alpha}}{(|\zeta|-1)^{1+\varepsilon}}}d \xi d \eta \leq C R,$$

\noindent for $\varepsilon < 2 \alpha$.
\QED
\end{pf}

\section{Proof of Theorem 2.}

Before proceeding to the proof we need to mention a result on quasiconformal mappings preserving Carleson measures. We say that $F:\mathbb{R}^2_{+} \rightarrow \mathbb{R}^2_{+}$ preserves Carleson measures if given any Carleson measure $\mu$ in $\mathbb{R}^2_{+}$, the measure $\nu$ defined in $\mathbb{R}^2_{+}$ as

\begin{equation}
   \nu(E)=\int_{F^{-1}(E)}a_F(z)d\mu(z) \nonumber
\end{equation}
\noindent is a Carleson measure, where

\begin{equation}
  a_F(z)=\displaystyle{\frac{1}{|B_z|}}\int_{B_z} J_F(\zeta)^{1/2}d\xi d\eta, \nonumber
\end{equation}
\noindent $B_z=B_z(1/2\;y)$ and $J_F$ is the Jacobian of $F$. The function $a_F(z)$ is somehow the quasiconformal substitute of $|F'|$ in Koebe distortion theorem. In fact $\operatorname{Im} F(z)\simeq a_F(z) y$ \cite{astala1}.\\

If we consider a quasiconformal map $\rho$ from $\mathbb{R}^2_{+}$ onto itself, it is well known that $\rho$ preserves Carleson measures if and only if $\rho|_{\mathbb{R}}$ is strongly quasisymmetric, i.e., it is locally absolutely continuous and $|\rho'|_{\mathbb{R}}| \in A_{\infty}$ \cite{astala}, \cite{mjose}. \\

We will need an analogous result for a quasiconformal map $\rho$ from $\mathbb{R}^2_{+}$ onto a domain bounded by a chord-arc curve.
%Since we could not find it in the literature, we state it here as a lemma.

\begin{lem}
\label{lema}
Let $\rho: \mathbb{R}^2_{+} \rightarrow \Omega$ a quasiconformal map and let $\Gamma=\rho(\mathbb{R})$ be a chord-arc curve. If $\rho|_{\mathbb{R}}$ is locally absolutely continuous and $|\rho'|_{\mathbb{R}}| \in A_{\infty}$, then $\rho$ preserves Carleson measures.
\end{lem}

%The proof would follow the same steps as in  $\Omega=\mathbb{R}^2_{+}$ and therefore we omit the details. In our case, a sequence of points $\{z_n\} \subset \Omega$ is interpolating  for $H^\infty (\Omega)$ if and only if
%
%\begin{enumerate}
%  \item $\exists \; \alpha >0, \; \; \forall j \neq k \;\; |z_j -z_k|\geq \alpha \; \mbox{dist}(z_j,\partial \Omega)$,
%  \item $\sum_{j\geq 0} \mbox{dist}(z_j,\partial \Omega)\;  \delta_{z_j} $ is a Carleson measure related to $\partial \Omega$,
%\end{enumerate}
%\noindent  where $\mbox{dist}(z,\partial \Omega)$ is the distance from the point $z$ to the curve $\partial \Omega$ and $\delta_z$ denotes the Dirac delta funcion \cite{zinsmeister}.\\
\begin{pf}
Note that since $\Gamma$ is chord-arc, there exists a global bilipschitz map $h$ that sends $\mathbb{R}^2_+$ onto $\Omega$ (\cite{pommerenke}, Theorem 7.9). As $h$ is bilipschitz, so is $h^{-1}$. Then $h^{-1} \circ \rho$ is absolutely continuous and $|(h^{-1}\circ \rho)'|_{\mathbb{R}}|\in A_{\infty}$. Since $h^{-1}\circ \rho$ sends $\mathbb{R}^{2}_+$ onto itself, the lemma follows from the previous case, that is $\Omega=\mathbb{R}^2_+$.
\end{pf}
\QED
\\

Let us now estate more precisely the result mentioned at the introduction due to Semmes (\cite{SEMMES2}, Theorem 0.1) and MacManus (\cite{macmanus}, Theorem 6.3): Let $\rho:\mathbb{C}\rightarrow\mathbb{C}$ be a quasiconformal mapping with dilatation $\mu$. Set $\tau=|\mu|^2/|y|$. If there exists $\gamma_0$ such that if $\|\tau\|_C\leq \gamma_0$ then $\rho(\mathbb{R})$ is a chord-arc curve, $\rho|_{\mathbb{R}}$ is absolutely continuous and $\log \rho' \in \operatorname{BMO}(\mathbb{R})$, with $\|\log \rho'\|_* \leq c \|\tau\|^{1/2}_C$ for some constant $c>0$.\\

We will apply this result when the measure $|\mu|^2/|y|$ is a vanishing Carleson measure. In fact, in this case not only $\log \rho' \in \operatorname{BMO}(\mathbb{R})$ but also $|\rho'|_{\mathbb{R}}| \in A_\infty$ as the following lemma shows.\\

\begin{lem}
\label{pesoinfinito}
Let $\rho:\mathbb{C}\rightarrow \mathbb{C}$ be a quasiconformal mapping with complex dilatation $\mu$. If $\mu$ has compact support and $|\mu|^2/|y|$ is a vanishing Carleson measure relative to $\mathbb{R}$, then $\rho(\mathbb{R})$ is a chord-arc curve, $\rho|_{\mathbb{R}}$ is absolutely continuous and $|\rho'|_{\mathbb{R}}| \in A_\infty$.
\end{lem}
\begin{pf}
To prove this, we may assume that $\mu=0$ outside a band $|y|<\varepsilon$ for some $\varepsilon >0$. Indeed, let $\tilde{\rho}$ be the solution of the Beltrami equation $\tilde{\rho}_{\bar{z}}/\tilde{\rho_{z}}=\mu(z)$ for $|y|<\varepsilon$ and $\tilde{\rho}_{\bar{z}}/\tilde{\rho_{z}}=0$ otherwise. Then $\rho=F \circ \tilde{\rho}$ where $F$ is conformal in the quasidisc $\tilde{\rho}(\{z:|y|<\varepsilon\})$.  We can then replace $\rho$ by $\tilde{\rho}$ in the whole plane without any loss of generality. Next, note that by choosing $\varepsilon>0$ small enough, we can make the Carleson norm of the measures as small as we want and, by the properties of $A_\infty$ weights mentioned in the preliminaries, we can conclude that $|\rho'|_{\mathbb{R}}| \in A_\infty$.\\
\QED
\end{pf}

To prove Theorem 2 we will follow Semmes approach in \cite{SEMMES}. Let $\Gamma$ be locally rectifiable quasicircle in the plane. Given a function $g$ defined on $\Gamma$, consider its Cauchy integral

\begin{equation}
   C_{\Gamma}(g)(z)=\displaystyle{\frac{1}{2 \pi i}}\int_{\Gamma} \displaystyle{\frac{g(\omega)}{\omega - z}} d \omega, \; \; z \notin \Gamma. \nonumber
\end{equation}

 We define the jump of $G=C_{\Gamma}(g)$ across $\Gamma$ at the point $z$ as $g_{+}-g_{-}$, where $g_{+}$ and $g_{-}$ denote the boundary values of $G$. As the classical Plemelj formula states,

\begin{equation}
  g_{\pm}(z)=\pm \displaystyle{\frac{1}{2}}g(z)+\displaystyle{\frac{1}{2 \pi i}}P.V. \int_{\Gamma} \displaystyle{\frac{g(\omega)}{\omega-z}}d\omega, \;\; z \in \Gamma. \nonumber
\end{equation}

\noindent Hence $g_{+}(z)-g_{-}(z)=g(z)$. Also, $G$ is holomorphic off $\Gamma$, so that $\bar{\partial}G=0$ on $\mathbb{C} \backslash \Gamma$.\\

Applying Green's theorem, we can reexpress these two conditions by \mbox{saying} that, in the distributional sense, $\bar{\partial} G=g d z_{\Gamma}$ on $\mathbb{C}$. This means that for any $\varphi \in C^{\infty}(\mathbb{C})$ with compact support

\begin{equation}
  \int_{\mathbb{C}} G(z)\bar{\partial} \varphi (z) dz \wedge d \bar{z}=-\int_{\Gamma} \varphi(z)(g^{+}(z)-g^{-}(z))dz_{\Gamma}, \nonumber
\end{equation}

\noindent where $dz_{\Gamma}$ denotes the usual measure on $\Gamma$ and $dz \wedge d\bar{z}$ represents the wedge product which is equal to $2$i$dxdy$.\\

%We can reexpress these two conditions by saying that, in the distributional sense, $\bar{\partial}G=g\;dz_{\Gamma}$ on $\mathbb{C}$, that is
%
%\begin{equation}
%  \int_{\Gamma}g\;dz_{\Gamma}=\int_{\Omega}(\bar{\partial}G)\;d\bar{z}\wedge dz,\nonumber
%\end{equation}
%
%\noindent where $dz_{\Gamma}$ is the usual measure on $\Gamma$ and $\partial \Omega=\Gamma$.\\

We can also say that $G$ is determined by the equation $\bar{\partial}G=0$ on $\mathbb{C} \backslash \Gamma$ and the condition $\mbox{jump}(G)=g$ on $\Gamma$, as if $\tilde{G}$ were another function with the same properties, then $\bar{\partial}(G-\tilde{G})=0$ in the sense of distributions and therefore, by Weyl's lemma, $G-\tilde{G}$ would be entire and a mild condition at $\infty$ would force it to be $0$.\\

Let $\tilde{G}=G \circ \rho$ on $\mathbb{C} \backslash \mathbb{R}$ and $f=g \circ \rho$ on $\mathbb{R}$, where $\rho:\mathbb{C} \rightarrow \mathbb{C}$ is a quasiconformal map that takes $\mathbb{R}$ into $\Gamma$. Then, $\bar{\partial}G=0$ off $\Gamma$ transforms into $(\bar{\partial}-\mu\partial)\tilde{G}=0$ off $\mathbb{R}$ with $\mbox{jump}(\tilde{G})=f$ across $\mathbb{R}$, where $\mu=\mu_\rho$. Again, in the distributional sense, we can say that $(\bar{\partial}-\mu\partial)\tilde{G}=f dx$.\\

 In order to prove the following result it  is convenient to change the pro\-blem a bit more. Let $F=C_{\mathbb{R}}(f)$ the Cauchy integral on $\mathbb{R}$ of f. Thus, $F$ is holomorphic off $\mathbb{R}$ and its jump across $\mathbb{R}$ is given by $f$, i.e., $\bar{\partial}F=fdx$.\\

 Let us now define $H=\tilde{G}-F$ on $\mathbb{C}
\backslash \mathbb{R}$. Then, $H$ has no jump across $\mathbb{R}$
and $\bar{\partial}H=\mu\partial \tilde{G}$. As $H$ has no jump, we can
consider that the previous equation holds on all of $\mathbb{C}$
in the sense of distributions (that there is no boundary piece).

\begin{thmB}
 \label{theoremB}
 Let $\rho$ be a quasiconformal map of the plane onto itself whose complex dilatation $\mu$ has compact support and satisfies that for some $\varepsilon> 0$ , \newline $|\mu|^2/|y|^{1+\varepsilon}$ is a Carleson measure relative to $\mathbb{R}$. Let $\Omega_{+}$ and $\Omega_{-}$ denote the two regions bounded by the quasicircle $\Gamma=\rho(\mathbb{R})$ and let $g \in L^{\infty}(\Gamma)$. Then $C_{\Gamma}(g) \in H^{\infty}(\Omega_{\pm})$ if and only if $C_{\mathbb{R}}(f) \in H^{\infty}(\mathbb{R}^2_{\pm})$ respectively, where \mbox{$f=g\circ \rho$.}
 \end{thmB}

%\begin{thmB}
%  Let $\Gamma$ be a $C^{1+\alpha}$ curve and let $g \in L^{\infty}(\Gamma)$. Then, the Cauchy integral $C_{\Gamma}(g) \in L^{\infty}(\mathbb{C})$ if and only if $C_{\mathbb{R}}(f) \in L^{\infty}(\mathbb{C})$, where $f$ denotes the pullback of $g$ under the conformal mapping.
%\end{thmB}

\begin{pf}
Set $G=C_{\Gamma}(g)$ where $g \in L^{\infty}(\Gamma)$ and $\tilde{G}=G\circ \rho \in L^{\infty}(\mathbb{R}^2_+)$. To show that $F=C_{\mathbb{R}}(f) \in H^{\infty}(\mathbb{R}^{2}_{+})$ we need to prove that $F|_{\mathbb{R}}\in L^{\infty}(\mathbb{R})$. Using the notation above, and since $H=\tilde{G}-F$, this is equivalent to prove that $H|_{\mathbb{R}} \in L^{\infty}(\mathbb{R})$. Since $\bar{\partial}H=\mu \partial \tilde{G}$, $\mu$ with compact support,
\begin{equation}
  H(a)=\displaystyle{\frac{1}{\pi i}}\int_{\mathbb{C}}
  \displaystyle{\frac{\bar{\partial}H(z)}{z-a}}dx dy=\displaystyle{\frac{1}{\pi i}}\int_{\mathbb{C}}\displaystyle{\frac{\mu(z)\partial \tilde{G}(z)}{z-a}}dx dy\;\; \text{for} \; \; a \in
  \mathbb{R}.
\end{equation}
\noindent For $a \in \mathbb{R}$ and $k$ an integer, let us denote $B_k=B_a(2^{-k})$. Then

\begin{eqnarray}
\label{acotacion2}
 |H(a)| & \lesssim & \sum_{k}2^{k+1}\int_{B_{k}\backslash B_{k+1}}|\mu(z)||\partial \tilde{G}(z)|dx dy \nonumber \\
         & \lesssim &  \sum_{k}2^{k+1}\left(\int_{B_{k}}\displaystyle{\frac{|\mu(z)|^2}{|y|}}dx dy\right)^{1/2}\left(\int_{B_{k}}|\partial \tilde{G}(z)|^2|y| dx dy \right)^{1/2}
\end{eqnarray}

\noindent Since $\nu(z)=|\mu(z)|^2/|y|^{1+\varepsilon}$ is a Carleson measure,

\begin{equation}
\label{acotacion3}
   \int_{B_k}\displaystyle{\frac{|\mu(z)|^2}{|y|}}dx dy=\int_{B_k}\displaystyle{\frac{|\mu(z)|^2}{|y|^{1+\varepsilon}}}|y|^{\varepsilon}dx dy
     \leq \|\nu\|_{C}2^{-k(1+\varepsilon)}.
\end{equation}

If $\tau=|\mu|^2/|y|$, the above inequality shows that $\|\tau\|_C \leq \|\nu\|_C 2^{-k\varepsilon}$ and, therefore, that it is a vanishing Carleson measure. By lemma \ref{pesoinfinito}, the quasicircle $\Gamma=\rho(\mathbb{R})$ is actually a chord-arc curve with small constant, $\rho|_{\mathbb{R}}$ is absolute continuous and $|\rho'|_{\mathbb{R}}| \in A_{\infty}$.\\

On the other hand, since $G \in H^{\infty}(\Omega_+)$, then $|G'(\zeta)|^2 \delta_\Gamma(\zeta) d\xi d\eta$ is a Carleson measure re\-la\-tive to $\Gamma$ with norm $\leq C \|g\|_{*}^2$ (\cite{SEMMES},Theorem 5.1). \\

%Let $Q$ be a Carleson cube in $\mathbb{R}^2_+$.
\noindent By lemma \ref{lema},
%{\bf Since $\rho|_{\mathbb{R}}$ is absolutely continuous and }$|\rho'|\in A_\infty$, by lemma \ref{lema}

\begin{equation}
   \sigma (B_k)=\int_{\rho(B_k)} a_{\rho^{-1}}(\zeta)|G'(\zeta)|^2\delta_{\Gamma}(\zeta)d\xi d\eta \lesssim 2^{-k}\nonumber
\end{equation}

\noindent So, the measure $\sigma$ defined in $\mathbb{R}^{2}_{+}$ is also a Carleson measure. By the Koebe distortion theorem for quasiconformal mappings \cite{astala1}, $a_{\rho^{-1}}(\zeta)\delta_{\Gamma}(\zeta) \simeq |\text{Im}(\rho^{-1}(\zeta))|=|y|$. Then

\begin{eqnarray}
\int_{B_k}|\partial \tilde{G}(z)|^2\;|y|\;dxdy &=& \int_{B_k} |G'(\rho(z))|^2|\partial \rho (z)|^2\;|y|\;dxdy \simeq \nonumber\\
&\simeq&\int_{B_k}|G'(\rho(z))|^2J_{\rho}(z)\;|y|\; dxdy \simeq\\
&\simeq &\int_{\rho(B_k)}|G'(\zeta)|^2a_{\rho^{-1}}(\zeta)\delta_{\Gamma}(\zeta)d\xi d\eta \lesssim 2^{-k}.\nonumber
\end{eqnarray}

%\begin{equation}
%  \sigma(B_k)\simeq \int_{B_k} |\partial \tilde{G}|^2\;|y|\;dxdy. \nonumber
%\end{equation}

\noindent This shows that $\lambda=|\partial \tilde{G}|^2 |y|\;dxdy$ is a Carleson measure relative to $\mathbb{R}$ and
\begin{equation}
\label{acotacion4}
\int_{B_{k}}|\partial \tilde{G}(z)|^2\;|y|\; dx dy \leq \|\lambda\|_C\; 2^{-k}.
\end{equation}

By (\ref{acotacion2}), (\ref{acotacion3}) and (\ref{acotacion4}) $|H(a)|\leq C(\|\nu\|_C,\|\lambda\|_C)\sum_{k}\left(2^{-\varepsilon/2}\right)^k <  \infty $ as we wanted to prove.\\

 Conversely, if $F$ were bounded, the same
argument would show that $\tilde{G}$ is bounded on $\mathbb{R}$ and then,
that $G \in H^{\infty}(\Omega_+)$.
\QED
\end{pf}

\noindent To conclude this paper we would like to propose the following problem:

\begin{problem}
Find conditions on $\mu$ so that Theorem 2 holds.
\end{problem}

%\bibliographystyle{model1a-num-names}
%\bibliography{
%\noindent {\bf References}

%\cite[p\'ag. 80]{Gratzer}
%\displaystyle{\lim_{ x \rightarrow 0}} f(x)
\end{document}